\documentclass{amsart}

\newtheorem{theorem}{Theorem}[section]
\newtheorem{lemma}[theorem]{Lemma}
\newtheorem{proposition}[theorem]{Proposition}
\theoremstyle{definition}
\newtheorem{definition}[theorem]{Definition}
\newtheorem{example}[theorem]{Example}
\theoremstyle{remark}
\newtheorem{remark}[theorem]{Remark}

\numberwithin{equation}{section}
\begin{document}

%

%
%

\title{Logarithmic singularities of solutions \\
to nonlinear partial differential equations  }

\author{Hidetoshi  TAHARA}
\address{Department of mathematics, Sophia University, 
  Kioicho, Chiyoda-ku, Tokyo 102-8554, Japan}  
\email{h-tahara@hoffman.cc.sophia.ac.jp}

\author{Hideshi YAMANE}
\address{Department of Physics,  Kwansei Gakuin 
University, Gakuen 2-1, Sanda, Hyougo 669-1337, Japan}
\email{yamane@ksc.kwansei.ac.jp}
    \thanks{
  This research was partially supported by Grant-in-Aid for 
  Scientific Research (No.16540169, No.17540182), Japan 
  Society for the Promotion of Science. 
  Parts of this work has been done during the authors' 
  stay at Wuhan University. They thank Professor 
  Chen Hua for hospitality and fruitful discussions.
}




\begin{abstract}
We construct a family of singular solutions 
to some nonlinear partial differential equations 
which  have resonances in the sense of 
a paper due to T.~Kobayashi. 
The leading term of a solution in our family 
contains a logarithm, possibly  multiplied by a monomial. 
As an application, we study nonlinear wave equations 
with quadratic nonlinearities. 
The proof is by the reduction to a Fuchsian equation with 
singular coefficients. 
\end{abstract}

\maketitle
\markboth{H.TAHARA and H.YAMANE}{LOGARITHMIC SINGULARITIES}

\part*{Introduction}

In this paper, we study singular solutions to 
nonlinear partial differential equations  with 
holomorphic (or real-analytic) coefficients. 
The  solutions to be constructed shall be 
 singular along  a \textit{non}characteristic hypersurface. 
This phenomenon presents a striking contrast to linear theory. 
Probably, the most   well-known example in this direction 
is  the KdV equation: 
\begin{equation}
\label{eq:KdV}
  u_{ttt}-6uu_x +u_x=0
  \quad(t, x\in \mathbb{C}).
\end{equation}
The surface $t=0$ is noncharacteristic but (\ref{eq:KdV}) 
has solutions of the form 
\[
  u=\frac{2}{t^2}
  +gt^2  +  ht^4 -\frac{1}{24}g_x t^5+\dots, 
\]
where $g=g(x)$ and $h=h(x)$ are arbitrary  functions. 
Note that many solutions have been obtained in the 
form of Laurent series for some integrable PDEs 
(a useful reference is  \cite{AC}).

In \cite{KS}, Kichenassamy and Srinivasan introduced 
an expansion of a generalized form in order to solve 
PDEs with polynomial nonlinearities. 
In their paper, 
the solutions  behave asymptotically 
\[
  u(t, x)\sim u_0 (x) t^\nu
  \quad (\text{as } t\to 0), 
\]
where $\nu$ is a rational number. The remainder term may 
contain logarithms. Besides this general result, 
specific cases are dealt with in  \cite{KL1}, \cite{KL2} and 
\cite{KichGowdy}: in these papers,   
it is proved that the Liouville equation $\square u=e^u$ and 
Einstein's vacuum equations  admit solutions led by 
logarithmic terms. 
 
On the other hand, 
in \cite{Kobayashi},  Kobayashi considered a certain kind of 
nonlinear PDEs, mainly those with polynomial nonlinearities, 
and constructed solutions of the form 
\[
  u(t, x)=t^{\sigma_c}\sum_{k=0}^\infty u_k (x)t^{k/p}, 
\]
where $\sigma_c\in\mathbb{Q}, p\in\mathbb{N}^*=\{1, 2, \dots\}$.  
The exponent $\sigma_c$, which is 
called the \textit{characteristic exponent},  is determined by 
 the nonlinear term of the equation and Kobayashi imposed a kind of 
generalized nonresonance condition on it. In particular, 
he assumed  $\sigma_c \ne 0, 1, 2, \dots, m-2$, where 
$m$ is the order of the equation. 
In the present paper, the authors shall deal with 
these excluded cases and construct solutions with 
a logarithm in the leading term. 
If $\sigma_c=l\in\{0, 1, 2, \dots, m-2\}$, then 
the asymptotic behavior of the solutions is 
\begin{equation}
\label{eq: u_small_o}
  u(t, x)\sim a(x)t^l \log t
   \quad (\text{as } t\to 0) 
\end{equation}
and the remainder term involves an arbitrary holomorphic  
function in $x$. 
Note that the case where $\sigma_c=0$ has already been treated 
in \cite{Yam} in a different formulation. 
This result about nonlinear wave equations shall be 
improved in Part~\ref{part:Yamane}, \S\ref{sec:wave}. 

Note that Tahara extended Kobayashi's result  for 
first-order equations with entire nonlinearities in 
\cite{Tahara1} and  
\cite{Tahara2}. 

All the above mentioned   authors employ 
the method of Fuchsian Reduction: 
the leading terms  can be found   by  formal calculation 
 and the remainder terms are obtained by solving 
nonlinear  Fuchsian equations. 
Here the word "Fuchsian" contains some ambiguity, 
because there are many versions of the notions of 
Fuchsian  or related equations. 
One has to choose a suitable version on each occasion.  
In the present work, we employ still another version, 
i.e.  \textit{equations with singular coefficients}.

The organization of the present paper is as follows: 
In Part~\ref{part:Yamane}, we shall study an 
equation  of the form 
\[
  \partial_t^m u =f
   \left(
   	t, x, (\partial_t^j \partial_x^\alpha u)
   \right), 
\]
where $f$ is holomorphic (real-analytic) in its arguments, 
and construct solutions with the asymptotic behavior 
$ u(t, x)\sim a(x)t^l \log t, \, l\in\{0, 1, \dots, m-2\}$. 
We shall explain 
how this equation is reduced to a 
Fuchsian equation with singular coefficients:  
\[
    ( t \partial_t  )^m u
      = F \Bigl( t,x,  (
                           ( t \partial_t  )^j
                           \partial_x^\alpha u
                       )
          \Bigr).
\]
Here $F(t, x, Z)$ is singular at $t=0$. 
The latter equation shall be solved in  Part~\ref{part:Tahara} 
by using the techniques developed in  \cite{yamazawa}. 
The characteristic exponents can be arbitrary.

\part{Logarithmic singularities}
\label{part:Yamane}

\section{Main result}
\label{sec:main}

 Let $(t,x)=(t,x_1,\ldots,x_n) \in \mathbb{C} \times \mathbb{C}^n$,
let $m \in \mathbb{N}^*$ be fixed and set:
$I_{m}=\{(j,\alpha) \in \mathbb{N} \times \mathbb{N}^{n} \,;\,
        j+|\alpha| \leq m \,\, \mbox{and} \,\, j<m \}$,
$N=\mbox{the cardinal of} \,\, I_{m}$, and
$U=(U_{j,\alpha})_{(j,\alpha) \in I_{m}} \in \mathbb{C}^N$. 
We set $\partial_t=\partial/\partial t$, 
\(
  \partial_x^\alpha=(\partial/\partial x_1)^{\alpha_1}
                    \dots
                    (\partial/\partial x_n)^{\alpha_n}  
\) for $\alpha=(\alpha_1, \dots, \alpha_n)$.

We study nonlinear PDEs of the form
\begin{equation}
\label{eq:main}
  \partial_t^m u
  =f
   \left(
   	t, x, (\partial_t^j \partial_x^\alpha u)_{(j, \alpha)\in I_m}
   \right)
   .
\end{equation}
Here $f(t, x, U)$ is holomorphic in 
$\{(t, x)\in\mathbb{C}_t\times\mathbb{C}^n_x; |t|< r_0, |x|< R_0\} 
\times \mathbb{C}^N_U$, 
where 
$r_0$ and $R_0$ are   positive constants. 
Note that Kobayashi assumed that 
$f$ was a polynomial in $U$. We shall deal with 
several infinite sums closely related to $f$. 
Their convergence follows from the fact that 
$f(t, x, U)$ is entire in $U$ 
(see Proposition~\ref{prop:a(x)} for example). 

We may write
\[
  f(t, x, U)=\sum_{\mu\in\mathcal{M}}
   f_\mu (t, x)U^\mu, 
   \quad
   \mu=(\mu_{j,  \alpha})_{(j, \alpha)\in I_m}, 
   \quad
   U^\mu=\prod_{(j, \alpha)\in I_m}U_{j, \alpha}^{\mu_{j,  \alpha}}
\]
for some subset $\mathcal{M}$ of $\mathbb{N}^N$, 
the set of $\mathbb{N}$-valued functions on $I_m$. 
We  assume that $f_{\mu}(t, x)$ does not vanish identically 
if $\mu\in\mathcal{M}$ (then  $\mathcal{M}$ is unique). 

We expand $f_\mu (t, x)$ in $t$:
\[
  f_\mu (t, x)=t^{k_\mu} \sum_{k=0}^\infty f_{\mu, k}(x)t^k.
\]
We assume that 
$f_{\mu, 0}(x)$ does not vanish identically. 
Summing up, we have 
\begin{equation}
\label{eq:expansion_of_f}
  f(t, x, U)
  =\sum_{\mu\in\mathcal{M}} f_\mu (t, x)U^\mu
  =\sum_{\mu\in\mathcal{M}}
  \left(
    t^{k_\mu} \sum_{k=0}^\infty
    f_{\mu, k}(x)t^k
  \right)U^\mu
  .
\end{equation}

We set 
\[
  |\mu|=\sum_{(j, \alpha)\in I_m} \mu_{j,  \alpha}, 
  \quad
  \gamma(\mu)=\sum_{(j, \alpha)\in I_m} j \mu_{j,  \alpha}
  .
\]

\begin{lemma}
\label{lem:gamma}
For any integer $l$, we have
\[ t^{\gamma(\mu)-l|\mu|}
 U^\mu
 =
  t^{\gamma(\mu)-l|\mu|}
  \prod_{(j, \alpha)\in I_m}
   U_{j, \alpha}^{\mu_{j,  \alpha}}
  =\prod_{(j, \alpha)\in I_m}
  (t^{j-l } U_{j, \alpha})^{\mu_{j,  \alpha}}
\]
\end{lemma}

\bigskip

We assume the following:
\begin{itemize}
\item[(A0)] 
\(\displaystyle \sup_{\mu\in \mathcal{M}, |\mu|\ge 2}
  \frac{\gamma(\mu)-m-k_\mu}{|\mu|-1}
  =l
  \in\{0, 1, 2, \dots, m-2\}. 
\)

Note that the left hand side is the characteristic exponent 
$\sigma_c$ in \cite{Kobayashi} 
(It is proved  in \cite{Kobayashi} that 
$\sigma_c<m-1$ holds if the supremum is attained by 
some $\mu$. )
\item[(A1)]
\(
  \mathcal{M}_0
  \underset{\text{def}}{=}
  \{\mu\in\mathcal{M}; 
     |\mu|\ge 2, \, 
     (\gamma(\mu)-m-k_\mu)/(|\mu|-1)=l
  \}
\)
is  non-empty: i.e.    the supremum in (A0) is attained. 
\item[(A2)] 
If $\mu\in\mathcal{M}_0$ and 
$\mu_{j\alpha}\ne 0$, then $j\ge l+1$  and $\alpha=0$. 
\item[(A3)] 
For a sufficiently small  positive constant $C>0$, 
we have 
\[
  m-l+k_\mu-\gamma(\mu)+l|\mu| 
  \ge C\sum_{ 
               \genfrac{}{}{0pt}{}{(j, \alpha)\in I_m}{j\le l}
            }\mu_{j, \alpha}
\] 
for any $\mu\in\mathcal{M}\setminus\mathcal{M}_0$. 
(This is trivial if $f$ is a polynomial by Lemma \ref{lem:m-l+k_mu} 
below. )
\end{itemize}

\begin{example}
\label{ex:proto}
The prototype is the ODE 
\[
  u^{(m)}= t^k \{u^{(m-1)}\}^2, \; 
  l=m-k-2.
\]
This equation is satisfied if 
$u^{(m-1)}=C_1  t^{-k-1}$, where $C_1$ is a suitable constant.  
In this case, we have $u\sim C_2 t^l \log t$ as $t\to 0$, 
where $C_2$ is another constant. 
\end{example}

We shall construct 
solutions to (\ref{eq:main}) 
which behave like (\ref{eq: u_small_o}). 
In order to give a precise statement, we 
introduce a function class $\widetilde{\mathcal O}_+$. 

We use the following notation:
\begin{itemize}
\item 
 ${\mathcal R}(\mathbb{C} \setminus \{0\})$,  the universal covering 
    space of $\mathbb{C} \setminus \{0\}$,
\item 
  $S_\theta = \left\{t \in {\mathcal R}(\mathbb{C} \setminus \{0\})
     \ ; \  |\mbox{arg}\,t| < \theta \right\}$,  a sector in 
    ${\mathcal R}(\mathbb{C} \setminus \{0\})$,
\item 
   $S(\varepsilon(y))  =  \left\{t \in 
   {\mathcal R}(\mathbb{C} \setminus \{0\});    
     0 < |t| < \varepsilon(\arg t) \right\}$,
     where $\varepsilon(y)$ is a positive-valued 
     continuous function on $\mathbb{R}_y$,
\item 
     $D_r = \left\{ x=(x_1, \ldots, x_n) \in  \mathbb{C}^n \ ; \  
      |x_i| < r \mbox{ for } i=1,\cdots,n \right\}$. 
\end{itemize}

\begin{definition}
$\widetilde {\mathcal O}_+$ denotes
the set of all $v(t,x)$ satisfying the following two conditions: 
\begin{itemize}
\item[$\mathrm{i})$]
	 $v(t,x)$ is a holomorphic 
	function on $S(\varepsilon (y)) \times D_r$ for some 
	positive-valued continuous function $\varepsilon (y)$ on $\mathbb{R}_y$ 
	and  $r>0$.
\item[$\mathrm{ii})$] 
	 there is an $a>0$  such that for any 
	$\tilde{r}\in]0, r[$ and $\theta >0$ we have
\[
     \max_{x \in D_{\tilde{r}}} |v(t,x)|= O(|t|^a) \quad 
     \mbox{(as $t \longrightarrow 0$ in $S_{\theta}$)}.
\]
\end{itemize}
\end{definition}

Our main result is the following:
\begin{theorem}
\label{thm:main}
Assume (A0)--(A3) and set $\beta_{j, l}=(-1)^{j-l-1} l! (j-l-1)!$ 
for $j\ge l+1$. 
Let  $A=a(x)$ be a solution to 
\begin{equation}
\label{eq:maintheorem_a(x)}
  \sum_{\mu\in\mathcal{M}_0} f_{\mu, 0}(x)  
  {\textstyle
   \left(
   	   \prod_{j=l+1}^{m-1}
  	  \beta_{j, l} ^{\mu_{j, 0}}
   \right)
   }
   A^{|\mu|-1}
       =\beta_{m, l}
    .
\end{equation}
Then,  for any holomorphic function $b(x)$ in a 
neighborhood of $x=0$,  
there exists a function $v(t, x)\in\widetilde{\mathcal O}_+$ 
such that 
\begin{align*}
   u(t, x)&=a(x)t^l \log t
   +t^l b(x)+t^l v(t, x)
   \\
   &=t^l
   \left\{a(x) \log t
   	+  b(x)+
   	 v(t, x)
    \right\}
\end{align*} 
is a solution to (\ref{eq:main}). 
\end{theorem}

\medskip

The convergence of the sum in the left hand side 
of (\ref{eq:maintheorem_a(x)}) 
shall be proved in Proposition~\ref{prop:a(x)}. 
Theorem~\ref{thm:main} itself shall be 
proved in \S\ref{sec:reduction}. 

We give some examples below. Note that 
the possibilities are 
$0\le l\le m-2$. 
See \S\ref{sec:wave} for an example of the case $m=2, l=0$. 

\begin{example}$[l=m-k-2;\ |\mu|=2, \ 
\gamma(\mu)=2(m-1)]$
\[
  \partial_t^m u=t^k (\partial_t^{m-1}u)^2
  +t^{k'} (\partial_t^j \partial_x^\alpha u)^2, 
\]
where 
\(
  j\le m-1,\ |\alpha|\le m-j, \
  2(m-1)-k\ge 2j-k'
\). This is just a PDE version of the 
ODE explained in Example~\ref{ex:proto}. 
The second term on the right hand side is 
a perturbation. Hence the set $\mathcal{M}_0$ consists 
of a single element and for the only $\mu\in\mathcal{M}_0$, 
we have $|\mu|=2, \ \gamma(\mu)=2(m-1)$. 
This is what is briefly stated after the semicolon 
between the square brackets. We shall employ the 
same shorthand notation in the following examples. 
\end{example}

\begin{example}$[m=3,\ l=0; 
\ \gamma(\mu)=3+k]$
\[
  \partial_t^3 u
  = t^k (\partial_t^2 u)^a 
  (\partial_t u)^{3+k-2a} u^b
  +t^{k'}
  (\partial_t^2 u)^p 
  (\partial_t u)^q (\partial_x^\alpha u)^r
  ,
\]
where $2p+q-k'-3<0, \  |\alpha|\le 3$. 
If the numerator $\gamma(\mu)-m-k_\mu$ vanishes 
in $\mathrm{(A0)}$, 
 the ratio also does  and  $|\mu|$ 
in the denominator is 
irrelevant. 
It makes it easy to construct   examples 
of the case  $l=0$. 
\end{example}

\begin{example}$[m=3,\ l=1; 
\ |\mu|=3, \ \gamma(\mu)=6]$
\[
  \partial_t^3 u
  = t(\partial_t^2 u)^3 
  +t^k (\partial_t^j \partial_x^\alpha u)^3, 
\]
where $3j-k\le 4,  \ |\alpha|\le 3-j$. 
\end{example}

\begin{example}$[m=4,\ l=0; 
\ \gamma(\mu)=k+4]$
\[
  \partial_t^4 
  = t^k (\partial_t^3 u)^a (\partial_t^2 u)^b (\partial_t u)^{4+k-3a-2b} 
  u^c
  +t^{k'} (\partial_t^3 u)^p (\partial_t^2 u)^q (\partial_t u)^r 
  (\partial_x^\alpha u)^s, 
\]
where $3p+2q+r-k'-4<0, \ |\alpha|\le 4$.
\end{example}

\begin{example}$[m=4,\ l=1; 
\ |\mu|=3, \ \gamma(\mu)=6]$
\[
  \partial_t^4 u
  = (\partial_t^3 u)^2 u
  + \partial_t^2 u \cdot \partial_t u \cdot \partial_x^\alpha u, 
\]
where $|\alpha|\le 4$. 
\end{example}

\begin{example}$[m=4, \ l=2; 
\ |\mu|=k+2, \ \gamma(\mu)=3(k+2)]$
\[
  \partial_t^4  u=t^k (\partial_t^3 u)^{k+2}
  +t^{k'} (\partial_t^j \partial_x^\alpha u)^{k+2}
  ,
\]
where $k'\ge k, \ 0\le j \le 2, \ |\alpha|\le 4-j$.
\end{example}

\begin{remark}
We constructed solutions 
with the growth order $|u|=O(|t^l \log t|)$ under the 
conditions (A0)--(A3). 
If $u$ is a solution with the growth order 
 $|u|=O(|t^{l'} \log t|)$, 
 $l'>l$, then we have 
  $|u|=O(|t^{(l+l')/2}|)$, 
 $(l+l')/2>l$ and 
 a result in \cite{Kobayashi} 
implies that  it can be extended as a 
holomorphic solution  up 
to some neighborhood of the origin.
\end{remark}

\section{Reduction to a Fuchsian equation}
\label{sec:reduction}
In this section, 
we reduce the equation  (\ref{eq:main}) to a 
Fuchsian equation with singular coefficients. 
The latter shall be the topic of 
Part~\ref{part:Tahara}. 

The assumption (A0)  is equivalent to the following:
\[
  m- l =\sup_{\mu\in \mathcal{M}, |\mu|\ge 2}
  (\gamma(\mu)- l  |\mu|-k_\mu) .
\]
Then we have   
\begin{equation}
\label{eq:sigma}
  m- l =\sup_{\mu\in \mathcal{M}}
  (\gamma(\mu)- l  |\mu|-k_\mu) 
  , 
\end{equation} 
because $\gamma(\mu)- l  |\mu|-k_\mu$ is smaller than 
$m- l $ if $|\mu|=0, 1$. 
Therefore, 

\begin{lemma}
\label{lem:m-l+k_mu}
We have
\begin{align*}
 &m- l +k_\mu= 
 \gamma(\mu)- l  |\mu| 
 \quad \text{if}\quad  \mu\in\mathcal{M}_0 , 
 \\
 &m- l + k_\mu>
 \gamma(\mu)- l  |\mu|  
 \quad \text{if}\quad \mu\in\mathcal{M}\setminus\mathcal{M}_0 .
\end{align*}
\end{lemma}

\begin{proposition} 
\label{prop:a(x)}
The sum in the left hand side of (\ref{eq:maintheorem_a(x)}) 
is an entire function in  $A\in\mathbb{C}$ and 
hence has at most one exceptional value in the 
sense of Picard's theorem in the value distribution 
theory of complex analysis. 
\end{proposition}

\begin{proof}
By  Lemma~\ref{lem:m-l+k_mu},   we have 
\begin{equation}
   \sum_{\mu\in\mathcal{M}_0} 
  f_{\mu, 0}(x) 
  \prod_{j=l+1}^{m-1}
   W_{j, 0} ^{\mu_{j, 0}}
 =
  r^{m-l}
  \sum_{\mu\in\mathcal{M}_0} 
  f_{\mu, 0}(x) 
  r^{k_\mu}
  \prod_{j=l+1}^{m-1}
  (W_{j, 0}/r^{j-l})^{\mu_{j, 0}}
\end{equation}
for any $W=(W_{j, \alpha})_{(j, \alpha)}$ and  $r>0$. 
Except for the factor $r^{m-l}$, the sum in the right hand side 
is nothing but a partial sum of (\ref{eq:expansion_of_f}), 
evaluated at $t=r, U=W/r^{j-l}$ by (A2). 
Hence it is convergent  
if $r>0$ is sufficiently small and 
so is the left hand side. 
Set  $W_{j, 0}=\beta_{j, l}A$, then we have
\[
  \sum_{\mu\in\mathcal{M}_0} 
  f_{\mu, 0}(x) 
  \prod_{j=l+1}^{m-1}
   W_{j, 0} ^{\mu_{j, 0}}
  =
  \sum_{\mu\in\mathcal{M}_0} f_{\mu, 0}(x)  
   \left(
   {\textstyle
	   \prod_{j=l+1}^{m-1}
  	  \beta_{j, l} ^{\mu_{j, 0}}
   }
   \right)
    A^{|\mu|}
    .
\]
This sum is convergent for any $A$. 
\end{proof}

\medskip

We set 
\[
	[\rho; j]
	=\Gamma(\rho+1)/\Gamma(\rho-j+1)
	=\rho(\rho-1)(\rho-2)\dots(\rho-j+1).
\]
Note that $[\rho; j]=0$ if $j-\rho$ is a positive integer. 
We define the sequence 
$\{b_{j, l}\}_j$ by 
\[b_{0, l}=0, \; b_{j+1, l}=[l; j]+(l-j)b_{j, l} .\]  
Then we have 

\begin{lemma}
\label{lem:devivatives_of_log}
\[
 \partial_t^j (t^l \log t)=
 \begin{cases}
t^{l-j}\{[l; j]\log t+ b_{j,l}\}
    & 
  (j\le l), 
  \\
  \beta_{j, l}t^{-(j-l)} 
  & (j\ge l+1).
 \end{cases}
 \]
\end{lemma}

We introduce  new unknown 
functions $a(x)$ and  $v(t, x)$ by setting 
\[
  u(t, x)=a(x)t^l \log t
  +t^l b(x)+t^l v(t, x). 
\]
Here the function $b(x)$ is arbitrary.

By using 
$t^{j-l}\partial_t^j t^l =[t\partial_t +l; j]$ 
and $t^l\in \mathrm{Ker\,}\partial_t^{l+1}$, we obtain

\begin{lemma}
\label{lem:u_to_v}
\begin{align*}
  t^{j-l} \partial_t^j \partial_x^\alpha u
  =
 \begin{cases}
  \partial_x^\alpha a (x)
  \{[l; j] \log t
              +b_{j, l}
            \}
            +  [t\partial_t +l; j] \partial_x^\alpha (b+v)
           &(j\le l),  
  \\
  \beta_{j, l} \partial_x^\alpha a (x)  
     +  [t\partial_t +l; j]  \partial_x^\alpha v
           &(j\ge l+1).
 \end{cases}
\end{align*}
\end{lemma}

Let us calculate the right hand side of (\ref{eq:main}). 
We have
\[
  f
   \Bigl(
   	t, x, (\partial_t^j \partial_x^\alpha u)_{(j, \alpha)\in I_m}
   \Bigr)
  =
   \sum_{\mu\in\mathcal{M}} 
    \left(
     \sum_{k=0}^\infty
    f_{\mu, k}(x)t^{k_\mu+k}
  \right)
   \prod_{(j, \alpha)\in I_m}
   (\partial_t^j \partial_x^\alpha u)^{\mu_{j,  \alpha}}
   .
\]

We extract the terms of the smallest weight and set 
\[
   S    =
   \sum_{\mu\in\mathcal{M}_0} f_{\mu, 0}(x) t^{k_\mu}
   \prod_{(j, \alpha)\in I_m}
   (\partial_t^j \partial_x^\alpha u)^{\mu_{j,  \alpha}}
   .
\]
Then by (A2), we have
\[
  S
  =
   \sum_{\mu\in\mathcal{M}_0} f_{\mu, 0}(x) t^{k_\mu}
   \prod_{j=l+1}^{m-1}
   (\partial_t^j   u)^{\mu_{j, 0}}
\] 
and it is free of logarithms. 
It is  a partial sum of $f(t, x, U)$ and 
its convergence is obvious. 
Note that  
\[
	m-l+k_\mu=\gamma(\mu)-l |\mu|
	=\sum_{j=l+1}^{m-1}(j-l)\mu_{j0}
\]
 for $\mu\in\mathcal{M}_0$. 
Hence  by Lemma~\ref{lem:gamma} we have
\begin{align*}
  t^{m-l} S
  &=\sum_{\mu\in\mathcal{M}_0} f_{\mu, 0}(x)  
   \prod_{j=l+1}^{m-1}
   (t^{j-l} \partial_t^j   u)^{\mu_{j,  0}}
  \\
  & = \sum_{\mu\in\mathcal{M}_0}  f_{\mu, 0}(x)  
   \prod_{j=l+1}^{m-1}
   \bigl\{\beta_{j, l} a(x)+[t\partial_t +l; j]  v\bigr\}
   ^{\mu_{j,  0}}
   .
\end{align*}
This quantity consists of terms of weight $0$. 
All the remaining parts of $t^{m-l} f$ 
consists of terms of 
positive weight (the weight of $\log t$ is $0$). 

By binomial expansion, we obtain 
\[
  t^{m-l} S=T_0+T_1+T_2, 
\]
where 
\begin{align*}
  T_0
  &=
   \sum_{\mu\in\mathcal{M}_0} f_{\mu, 0}(x)  
   \prod_{j=l+1}^{m-1}
   \{\beta_{j, l} a(x)\}^{\mu_{j, 0}}
   =
   \sum_{\mu\in\mathcal{M}_0} f_{\mu, 0}(x)  
   \Biggl(
   {
	   \prod_{j=l+1}^{m-1}
  	  \beta_{j, l} ^{\mu_{j, 0}}
   }
   \Biggr)
    a(x)^{|\mu|}
    , 
  \\
  T_1
  &= 
   \sum_{\mu\in\mathcal{M}_0} f_{\mu, 0}(x)  
  \sum_{j=l+1}^{m-1} 
  \Biggl(
  {
	   \prod_{i\ne j}
	   \{\beta_{i, l} a(x)\}^{\mu_{i, 0}}
  }
  \Biggr)
  \mu_{j, 0}
    \{\beta_{j, l} a(x)\}^{\mu_{j, 0} -1}
   [t\partial_t +l; j] v
   \\
   &= 
   \sum_{\mu\in\mathcal{M}_0} f_{\mu, 0}(x)  
  \sum_{j=l+1}^{m-1} 
  \Biggl(
  { 
   	\prod_{i\ne j}
 	 \beta_{i, l}^{\mu_{i, 0}}
  }
  \Biggr)
  \mu_{j, 0}
    \beta_{j, l}^{\mu_{j, 0} -1}
    a(x)^{|\mu|-1}
   [t\partial_t +l; j] v
   , 
  \\
  T_2&=\text{a polynomial in  } (t\partial_t)^k v \;
  (k=0, 1, \dots,  m-1) 
  \text{ free of terms of degree} \le 1
  .
\end{align*}
Note that $T_0$ is free of $v$ and its derivatives. 

On the other hand, we have 
\[
  t^{m-l} \partial_t^m u 
  =\beta_{m, l}a(x)
  +    [t\partial_t+l; m]v
  .
\]

We multiply the left and right hand sides of 
(\ref{eq:main}) by $t^{m-l}$. 
If $a(x)$ is determined by (\ref{eq:maintheorem_a(x)}), 
the function 
$u=at^l \log t+t^l b +t^l v$ is a solution to (\ref{eq:main}) if 
and only if $v$ is a solution to the equation below:
\begin{equation}
\label{eq:v}
 [t\partial_t+l; m]v-T_1=t^{m-l}f-T_0-T_1
 .
\end{equation}

Set 
\[
  \delta(\mu)
  =m-l+k_\mu-\gamma(\mu)+l|\mu|, 
  \qquad
  |\mu|_l 
  =
  \sum_{
	       \genfrac{}{}{0pt}{}
	       {(j, \alpha)\in I_m}{j\le l}
	     }
	  \mu_{j, \alpha}
	  . 
\] 
Since $\delta(\mu)$ is an integer, 
Lemma~\ref{lem:m-l+k_mu} implies that 
 $\delta(\mu)=0$ for  $\mu\in\mathcal{M}_0$ 
and that $\delta(\mu)\ge 1$ 
for  $\mu\in\mathcal{M}\setminus\mathcal{M}_0$. 
Moreover,  (A3) can be written in a simple form:
\[
  \delta(\mu)
    \ge C |\mu|_l , 
    \qquad \mu\in\mathcal{M}\setminus\mathcal{M}_0 
    .
\]
These two estimates imply that    
\begin{equation}
\label{eq:delta_ineq}
  \delta(\mu)
  =
  \delta(\mu)/3+\delta(\mu)/3+\delta(\mu)/3
  \ge
  1/3+\delta(\mu)/3+C|\mu|_l/3 
\end{equation}
holds for any $\mu\in\mathcal{M}\setminus\mathcal{M}_0$. 

By the way, we have
\begin{align}
  \label{eq:t^(m-l) f}
  t^{m-l} f(t, x, U)
  &= t^{m-l}\sum_{\mu\in\mathcal{M}}
  \left(
    t^{k_\mu} 
    \sum_{k=0}^\infty f_{\mu, k}(x)t^k
  \right) 
  U^\mu
   \\
    \nonumber
   &=
  \sum_{\mu\in\mathcal{M}}
  \left(
       \sum_{k=0}^\infty f_{\mu, k}(x)
    t^{\delta(\mu)+ k}
         \right) 
   \prod_{(j, \alpha)\in I_m}
   (t^{j-l} U_{j, \alpha})^{\mu_{j,  \alpha}}
   .
\end{align}

\begin{proposition}
\label{prop:Tahara}
The equation (\ref{eq:v}) satisfies the conditions 
$\mathrm{C}_1)$ and $\mathrm{C}_2)$ in 
Part~\ref{part:Tahara}. 
Here our unknown function $v(t, x)$ plays the role of $u$ in 
Part~\ref{part:Tahara}. 
\end{proposition}

\begin{proof}
Set 
\begin{align*}
 &a^{(\alpha)} = \partial_x^{\alpha} a(x), \quad
  b^{(\alpha)} = \partial_x^{\alpha} b(x), 
\\
 &  W_{j,  \alpha}
    =t^{j-l}U_{j, \alpha}
    - a^{(\alpha)}  t^{j-l}\partial_t^j (t^l \log t)
    -[l; j]  b^{(\alpha)}, 
\end{align*}
then $W_{j,  \alpha}$  corresponds to 
\(
  t^{j-l}\partial_t^j \partial_x^\alpha (t^l v)= 
  [t\partial_t +l; j] \partial_x^\alpha v(t, x)
\). 
Note that  
$\{[t\partial_t +l; j] v\}_{j=0, \dots, m-1}$ is equivalent to 
$\{ (t\partial_t )^j v\}_{j=0, \dots, m-1}$, 
the latter being used in Part~\ref{part:Tahara}. 
They are transformed 
into each other by the action of a  lower triangular matrix 
whose diagonal elements are all $1$.

For brevity, we set 
\begin{align*}
  \widetilde{W}_{j, \alpha}
    &=W_{j,  \alpha}
    + a^{(\alpha)}   t^{j-l} \partial_t^j  (t^l \log t)
    +[l; j]  b^{(\alpha)}
    (=t^{j-l}U_{j, \alpha})
    \\
    &=
 \begin{cases}
 	W_{j,  \alpha}
    + a^{(\alpha)}   \{[l; j]\log t+ b_{j,l}\}
    +[l; j]  b^{(\alpha)}
    \quad
  (j\le l), 
    \\
    W_{j,  \alpha}
    + \beta_{j, l} a^{(\alpha)}  
        \quad
  (j\ge l+1).
 \end{cases}
\end{align*}
Here we have used Lemma~\ref{lem:devivatives_of_log} 
and the fact that $[l; j]=0$ if $j\ge l+1$. 
By (\ref{eq:t^(m-l) f}) we have 
\[
   \{t^{m-l} f(t, x, U)-T_0 -T_1\} -T_2
  =t^{m-l} f(t, x, U)-t^{m-l} S
   =I_1 + I_2, 
\]
where 
\begin{align*}
 I_1&= \sum_{\mu\in \mathcal{M}_0 }
  \sum_{k\ge 1} 
  f_{\mu, k} (x)t^k 
  \prod_{j=l+1 }^{m-1} 
    	\widetilde{W}_{j, 0}^{\mu_{j, 0}}
    	,
  \\
 I_2 &=
   \sum_{\mu\not\in \mathcal{M}_0 }
       \sum_{k=0}^\infty f_{\mu, k}(x)
      t^{\delta(\mu)+k}
  \prod_{  
            \genfrac{}{}{0pt}{}
              {(j, \alpha) \in I_m }{ j\le l}  
        }
         \widetilde{W}_{j, \alpha}^{\mu_{j,  \alpha}}
       \prod_{
               \genfrac{}{}{0pt}{}
               {(j, \alpha)\in I_m}{j\ge l+1}  }
       \widetilde{W}_{j, \alpha}^{\mu_{j,  \alpha}}
       .
\end{align*}
The convergence of $I_1$ can be proved by 
the method of Proposition~\ref{prop:a(x)}. 
If $|t|\leq r$, then the following estimates holds: 
\[
  |I_1| \le 
  |t| r^{m-l-1} \sum_{\mu\in\mathcal{M}_0} \sum_{k\ge 1} 
  |f_{\mu, k} (x)| 
  r^{k_\mu + k} 
  \prod_{j=l+1}^{m-1}
  \Bigl(
      |\widetilde W_{j, 0}|/r^{j-l}
  \Bigr)^{\mu_{j, 0}} 
  .
\]
We see that 
$I_1$ and its derivatives in $W$ are of 
order $O(|t|)$.

Next let us consider $I_2$. 
The trivial fact $\delta(\mu)>0, \mu\not\in \mathcal{M}_0$ 
helps, to be sure, but it is not good enough. 
If $j\le l$, 
the quantity 
$\widetilde W_{j, \alpha}$ contains a logarithm,  
whose unboundedness is the 
greatest obstacle. We overcome it by assuming (A3). 
The trick is the following fact: 
if $C>0$, then 
\(
 t^C \widetilde{W}_{j, \alpha}
\) 
is bounded as $t\to 0$. 

If  $|t|<r<1$, 
the inequality (\ref{eq:delta_ineq}) implies that 
\begin{align*}
  |t^{\delta(\mu)}|
  &\le
  |t|^{1/3}
  \times
  (r^{1/3})^{\delta(\mu)}
  \times
  (|t|^{C/3})^{|\mu|_l}
    \\
  &=
    |t|^{1/3}
   r^{(m-l+k_\mu)/3}
    \times
    \prod_{   
             \genfrac{}{}{0pt}{}
             {(j, \alpha)\in I_m }{ j\le l}
          }
    \left(
    	\frac{|t|^{C/3}}{r^{(j-l)/3}}
    \right)^{\mu_{j, \alpha}}
    \times
    \prod_{  
             \genfrac{}{}{0pt}{}
             {(j, \alpha)\in I_m}{j\ge l+1}
          }
    \left(
    	\frac{1}{ r^{(j-l)/3}}
    \right)^{\mu_{j, \alpha}}
    .
\end{align*}
Therefore 
\begin{align*}
  |I_2|
  &\le
  |t|^{1/3}
  \sum_{\mu\not\in\mathcal{M}_0} \sum_{k=0}^\infty 
  |f_{\mu, k}(x)|   (r^{1/3})^{m-l+k_\mu}
  |t|^k
  \\
  &\qquad\qquad
  \times
 \prod_{       
               \genfrac{}{}
               {0pt}{}{(j, \alpha)\in I_m}{j\le l}
       }
    \left(        r^{(l-j)/3}
    	         |t|^{C/3}  
    	          |\widetilde{W}_{j, \alpha} |
    \right)^{\mu_{j, \alpha}}
  \times
 \prod_{  
          \genfrac{}{}
          {0pt}{}{(j, \alpha)\in I_m}{j\ge l+1}
       }
    \left(
    	        r^{(l-j)/3}
    	        |\widetilde{W}_{j, \alpha}  |
    \right)^{\mu_{j, \alpha}}
    .
\end{align*}
If $r>0$ is sufficiently small, then 
$I_2$ is convergent in $|t|< r< r^{1/3}$. 
We see that 
$I_2$ and its derivatives in $W$ are of 
order $O(|t|^{1/3})$.
\end{proof}

\noindent
\textit{End of Proof of Theorem~\ref{thm:main}.} 
Theorem~\ref{thm:main} follows from 
Proposition~\ref{prop:Tahara} and 
Theorem~\ref{thm:existence_tahara} of Part~\ref{part:Tahara}. 

\begin{remark}
\label{rem: arbitrary}
There is an arbitrary function $b(x)$  in  
 the family of solutions in Theorem~\ref{thm:main}. 
 In some cases, 
there may be more: as is stated 
in Remark~\ref{remark:nonuniqueness} of  Part~\ref{part:Tahara}, 
the equation (\ref{eq:v}) may admit a family of solutions 
involving one or more arbitrary functions in $x$. 
\end{remark}

\section{Nonlinear wave equation}
\label{sec:wave}

We can relax the condition imposed in \cite{Yam}. 
Although we formulate our result in the complex domain, 
it is trivial that an analogous result holds in 
the real-analytic category. 

We consider 
\begin{equation}
\label{eq:nonlinear_wave}
 \square u(s, y)
 =g(s, y; u, \partial_s u, \nabla_y u)
\end{equation}
in an open set of $\mathbb{C}^{n+1}=\mathbb{C}_s \times \mathbb{C}^{n}_y$. 
Here 
\(
 \square=\partial^2/\partial s^2
 -\sum_{i=1}^n \partial^2/\partial y_i^2
\), \;
\(
  \nabla_y u 
  =(\partial u/\partial y_1, \dots, 
  \partial u/\partial y_n
  )
\). 
We assume that  $g(s, y; z, \sigma, \eta)$ is 
a holomorphic function in all its arguments and  
is entire in $(z, \sigma, \eta)$. 
Moreover we assume that it 
is a polynomial of degree $2$ in 
$(\sigma, \eta)$. Its homogeneous part 
of degree $2$ is 
denoted by $g_2$. 

Let $\psi(y)$ be a holomorphic function 
with 
\begin{equation}
\label{eq:Psi_ne_0}
   1-|\nabla_y \psi(y)|^2\ne 0
   , 
\end{equation}  
where
\( 
  \nabla_y \psi (y)=
  \bigl(\psi_1 (y), \dots, \psi_n (y)\bigr), 
  \;
   \psi_i (y)= \partial \psi(y)/\partial y_i  
  (i=1, 2,  \dots, n)
\). 
Moreover, assume that 
\begin{equation}
\label{eq:g_2_ne_0}
  g_2 \bigl(\psi(y), y; 0, 1,  -\nabla_y \psi (y)\bigr)\ne 0 
   .
\end{equation}
Note that this assumption corresponds to $k_\mu=0$, 
where $k_\mu$ is as in \S\ref{sec:main}. 

\begin{theorem}
\label{thm:nonlin_wave}
Assume 
(\ref{eq:Psi_ne_0}) and (\ref{eq:g_2_ne_0}). 
Then, in a neighborhood of the hypersurface $\{s=\psi(y)\}$, 
there  exists a family of solutions $u(s, y)$ to 
 (\ref{eq:nonlinear_wave}) with the 
 asymptotic behavior 
\[
   u(s, y)
   \sim 
   -\frac{1-|\nabla_y \psi(y)|^2}
         {g_2 \bigl(\psi(y), y
            ; 0, 1, -\nabla_y \psi (y)\bigr)
         }
   \log \bigl(s-\psi(y)\bigr). 
\]
Here the remainder term involves an arbitrary 
holomorphic function on $\{s=\psi(y)\}$. 
\end{theorem}

\begin{proof}
Set 
\[
  t=s-\psi(y), \quad x=y, \quad
  \Psi=1-|\nabla_y \psi(y)|^2 (\ne 0). 
\]
Then, as is proved in \cite{KL1}, we have
\begin{align*}
   & \partial_s=\partial_t, 
   \quad
   \partial_{y_i}=-\psi_i \partial_t + \partial_{x_i}, 
   \\
   &\square=\square_{s, y}
   =\Psi \partial_t^2 +
    2\sum_{i=1}^n \psi_i \partial_{x_i}\partial_t 
    +(\Delta_y \psi)\partial_t -\Delta_x
    .
\end{align*}
In a neighborhood of $t=s-\psi(y)=0$, 
the original equation (\ref{eq:nonlinear_wave}) 
becomes 
\[
  \partial_t^2 u
  =(\text{linear part})
  +\Psi^{-1} 
  g
  \Bigl(t+\psi(x), x; u, \partial_t u, 
     \bigl((-\psi_i \partial_t +\partial_{x_i}) u \bigr)
     _{i=1, \dots, n} 
  \Bigr)
  .
\]
When we expand the right hand side in 
a power series in $t$,  we find   the term 
\[
  \Psi^{-1}
  g_2 \bigl(\psi(x), x; 0, 1, -\nabla_y \psi(y)\bigr)
  (\partial_t u)^2
  .
\]
It corresponds to $\mu$ with 
$\mu_{1, 0}=2, \mu_{j, \alpha}=0 \;(\text{otherwise})$. 
For this $\mu$, we have 
\[k_\mu=0, \; \gamma(\mu)=2, \;
   f_{\mu, 0}(x)
 =\Psi^{-1}
 g_2
 \bigl(
   \psi(x), x; 0, 1, -\nabla_y \psi(y)
 \bigr)
 .
\]
It is the only element of $\mathcal{M}_0$ and 
we can apply Theorem~\ref{thm:main} with $l=0$. 
Since $\beta_{2, 0}=-1, \beta_{1, 0}=1$, the equation 
(\ref{eq:maintheorem_a(x)}) reduces to 
\(
  f_{\mu, 0}(x) A=-1
\). 

Theorem~\ref{thm:main}  enables us to construct a 
solution in a neighborhood of each point  on the hypersurface. 
In spite of 
Remark~\ref{rem: arbitrary}, these  solutions overlap, if 
$b(x)$ is fixed,  
because they are constructed in the same way,  i.e. by 
 Proposition~\ref{prop:a, 2a}, (\ref{eq:u_1}) and (\ref{eq:u_k}). 
\end{proof}

\part{Nonlinear Fuchsian equations with 
singular coefficients}
\label{part:Tahara}

We shall generalize the result of  \cite{yamazawa}
 to the case where the equations have singular 
coefficients at $t=0$. 
We employ the same notation as in Part~\ref{part:Yamane}. 
Two more sets have to be introduced: 
\begin{itemize}
\item $S_\theta (\delta) = \left\{t \in S_\theta \ ; \  
    0< |t| < \delta \right\}$ a sectorial domain in 
    ${\mathcal R}(\mathbb{C} \setminus \{0\})$,

\item $S_\theta (\varepsilon(y)) 
                        = S_\theta \cap S(\varepsilon(y))$.
\end{itemize}

\section{An existence theorem}
\label{sec:existence}

We consider 

\begin{equation}
\label{eq:nonlin_fuchs}
     ( t \partial_t  )^m u
      = F \Bigl( t,x, \bigl(
        ( t \partial_t )^j
         \partial_x^\alpha u
        \bigr)_{(j,\alpha) \in I_m}
          \Bigr)
\end{equation}
with the unknown function $u=u(t,x)$. Here the function 
$F$ is allowed to be singular at $t=0$. 
Typically, it may involve powers of $\log t$. 
More precisely, we assume:
\begin{itemize}
\item[$\mathrm{C}_1$)]
   $F(t,x,Z)$ is  a holomorphic function 
      in   $(t,x,Z)$, 
    $Z=(Z_{j,  \alpha})_{(j, \alpha)\in I_m}\in \mathbb{C}^N$ on 
   $S(\varepsilon (y)) \times D_{R_0} \times \{|Z|< L \}$ 
   for a positive-valued continuous 
   function $\varepsilon (y)$ and 
   constants $R_0>0$, $L>0$.
\item[$\mathrm{C}_2$)]  there exist a constant $s > 0$ 
    and holomorphic functions $c_j(x)$ ($0\le j\le m-1$) 
     on $D_{R_0}$ such that for any $\theta >0$, 
    $(j,\alpha) \in I_m$ and $(i,\beta) \in I_m$ we
    have
\begin{align*}
    \mathrm{i})& \enskip 
      \sup_{x \in D_{R_0}} |F(t,x,0)|=O(|t|^s) \enskip
      \mbox{(as $S_{\theta} \ni t \longrightarrow 0$)}, \\
    \mathrm{ii})& \enskip \sup_{x \in D_{R_0}} \Bigl| 
       \dfrac{\partial F}{\partial Z_{j,0}}
              (t,x,0) -c_j(x) \Bigr|=O(|t|^s)
       \enskip \mbox{(as $S_{\theta} \ni t
     \longrightarrow 0$)}, \\
    \mathrm{iii})& \enskip \sup_{x \in D_{R_0}} \Bigl| 
       \dfrac{\partial F}{\partial Z_{j,\alpha}}
              (t,x,0) \Bigr|=O(|t|^s)
       \enskip \mbox{(as $S_{\theta} \ni t
     \longrightarrow 0$)}  \enskip \mbox{if $|\alpha|>0$}, \\
    \mathrm{iv})& \enskip \sup_{x \in D_{R_0}, \,|Z| < L}\,
       \Bigl| \dfrac{\partial^2 F}
        {\partial Z_{j,\alpha} \partial Z_{i,\beta}}
              (t,x,Z) \Bigr|=O(1)
       \enskip \mbox{(as $S_{\theta} \ni t
     \longrightarrow 0$)}.
\end{align*}
\end{itemize}

\par
\smallskip
\noindent
   Then we have:
\begin{theorem}
\label{thm:existence_tahara}
     Assume the conditions $\mathrm{C}_1)$ and $\mathrm{C}_2)$. 
Then, the equation (\ref{eq:nonlin_fuchs}) has a solution $u(t,x)$  
in the class $\widetilde{\mathcal O}_+$.
\end{theorem}

\begin{remark}
\label{remark:nonuniqueness}
   The solution of  (\ref{eq:nonlin_fuchs}) 
in $\widetilde{\mathcal O}_+$ is not
necessarily unique.  There may be a family of 
solutions involving one or more arbitrary 
functions in $x$. See \cite{yamazawa}.
\end{remark}

Note that this theorem is essential in the proof of 
Theorem~\ref{thm:main}. 
\section{Some preparatory discussion}

   Before the proof of Theorem~\ref{thm:existence_tahara}, 
   let us present some preparatory discussion.     
For a function 
$\phi(x)$ on $D_r$,  we define the norm $\| \phi \|_r$ by
\[
     \| \phi \|_r = \sup_{x \in D_r}|\phi(x)|.
\]
   Let $\varepsilon (y)$ be a positive-valued continuous function
on $\mathbb{R}_y$.  We say that $\varepsilon (y)$ is decreasing in $|y|$
if the following condition holds: $|y_1| \leq |y_2|$ implies
$\varepsilon (y_1) \geq \varepsilon (y_2)$.

\begin{definition}
  (1) For $d \geq0$ and $\theta > 0$, we 
denote by $\widetilde{\mathcal O}_d(S_{\theta}(\varepsilon (y)) \times D_R)$ 
the set of all the holomorphic functions on 
$S_{\theta}(\varepsilon (y)) \times D_R$ that satisfy the
following estimate: for any $0<r<R$ there is a constant $C>0$ such that
\[
    |u(t,x)| \leq C |t|^d \quad
     \mbox{on} \quad  S_{\theta}(\varepsilon (y)) \times D_r.
\]
\noindent
   (2) We set \; 
\(  \widetilde {\mathcal O}_d(S(\varepsilon (y)) \times D_R)
   = \bigcap_{\theta>0} 
   \widetilde{\mathcal O}_d(S_{\theta}(\varepsilon (y)) \times D_R)
\). 
\end{definition}

\bigskip

   Let $m \in \mathbb{N}^*$ and $c_j(x)$ ($j=0,1,\ldots,m-1$) be as in 
   \S\ref{sec:existence}. 
Set
\begin{equation} 
\label{eq:char_poly} 
   C(\lambda,x) = \lambda^m -c_{m-1}(x) \lambda^{m-1} - \cdots-
                 c_{1}(x) \lambda -c_0(x),
\end{equation} 
and 
denote by $\lambda_1(x), \ldots, \lambda_m(x)$ the roots of 
$C(\lambda,x)=0$ in $\lambda$, and let us consider the following 
equation:
\begin{equation}
\label{eq:Cv=g}
    C \bigl( t \partial_t, x \bigr) v
    =g(t,x). 
\end{equation}

\begin{proposition}
\label{prop:a, 2a}
  Let $a>0$. Suppose that
\begin{equation}
\label{eq:g_estimate}
   \bigl\{ a, 2a, 3a, \ldots \bigr\} \cap 
   \bigl\{ \mathrm{Re} \lambda_1(0),\ldots, 
           \mathrm{Re} \lambda_m(0) \bigr\} 
   = \emptyset   
\end{equation}
and that $\varepsilon (y)$ is decreasing in $|y|$.
Then we can take a sufficiently small $R_1>0$ satisfying
the following properties $\mathrm{(*)}_k$ and $\mathrm{(\sharp)}_k$ for 
$k=1,2, \ldots:$.

\medskip\noindent
$\mathrm{(*)}_k$: 
   \, For any 
\(
  g(t,x) \in \widetilde {\mathcal O}_{ak}
  (S(\varepsilon (y)) \times D_{R_1})
\),  
the equation (\ref{eq:Cv=g}) has 
a solution $v(t,x) \in \widetilde {\mathcal O}_{ak}
(S(\varepsilon (y)) \times D_{R_1})$.

\medskip\noindent
$\mathrm{(\sharp)}_k$: \, Moreover, if $g(t,x)$ satisfies
\begin{equation}
\label{eq:v_estimate1}
    \| g(t, x) \|_r \leq C |t|^{ak} \quad 
    \mbox{on} \quad S_{\theta}(\varepsilon (y))
\end{equation}
for some $0<r<R_1$, $C>0$ and $\theta>0$, 
we have the estimate
\begin{equation}
\label{eq:v_estimate2}
   \|  (t \partial_t )^{j} 
            v(t) 
   \|_r 
   \leq 
    \frac{M_{\theta}}{k^{m-j}} \,C |t|^{ak} \enskip
    \mbox{on $S_{\theta}(\varepsilon (y))$}
    \enskip \mbox{for $j=0,1, \ldots, m-1$},
\end{equation}
where the constant $M_{\theta}>0$ is independent of $k$, 
$g(t,x)$, $r$ and $j$.
\end{proposition}
\begin{proof}
This  proposition   
can be proved by the same argument as in the proof of 
Lemma 6 in \cite{yamazawa}. 
We can choose a suitable path of integration because of 
the assumption that  $\varepsilon (y)$ is  
decreasing. 
\end{proof}

\begin{lemma}[Nagumo's Lemma]
\label{lem:nagumo}
   If $\phi(x)$ is a holomorphic function on 
$D_R$ and if 
\[
     \|\phi \|_r \leq \frac{C}{\,(R-r)^{\,b} \,}  \quad 
     \mbox{for any $0< r<R$}
\]
holds for some $C \geq 0$ and $b \geq 0$, then we have
\[
    \left\| \frac{\partial \phi}{\partial x_i} \right\|_r
    \leq \frac{e(b+1)C}{(R-r)^{\,b+1}}  
    \quad \mbox{for any \,$0<r<R$ \,
                         and \,$i=1, \ldots, n$}.
\]
\end{lemma}

\begin{proof}
See Nagumo \cite{Nagumo} or 
Lemma 5.1.3 of H\"{o}rmander \cite{Hormander}
\end{proof}

\section{Proof of Theorem~\ref{thm:existence_tahara}}

   Assume the conditions $\mathrm{C}_1)$ and $\mathrm{C}_2)$.
Then, by expanding $F(t,x,Z)$ in $Z$,  our equation 
(\ref{eq:nonlin_fuchs}) is
written in the form
\begin{align}
\label{eq:nonlin_fuchs2}
   C  ( t \partial_t, x  )u
   = &a(t,x)+ 
     \sum_{
             \genfrac{}{}{0pt}{}
             {(j,\alpha) \in I_m}{|\alpha|>0}
          } 
     b_{j,\alpha}(t,x)
     ( t \partial_t  )^j 
      \partial_x^{\alpha}u 
     \\
  &+\sum_{|\nu| \geq 2} g_{\nu}(t,x)
      \prod_{(j,\alpha) \in I_m} \left[
      (t \partial_t)^j 
        \partial_x ^{\alpha}u
     \right]^{\nu_{j,\alpha}},
     \nonumber
\end{align}
where $\nu=(\nu_{j,\alpha} )_{(j,\alpha) \in I_m} \in \mathbb{N}^N$
and $|\nu|=\sum_{(j,\alpha) \in I_m} \nu_{j,\alpha}\ge 2$. 
The coefficients  $a(t,x)$, $b_{j,\alpha}(t,x)$ 
and $g_{\nu}(t,x)$  are all holomorphic 
functions on $S(\varepsilon (y)) \times D_{R_0}$ 
with suitable growth order to be specified below. 
By replacing  $\varepsilon (y)$ 
if necessary, we may suppose that $0< \varepsilon (y) \leq 1$ 
and that $\varepsilon (y)$ is decreasing in $|y|$.

\bigskip

Let us construct a formal solution. 
By taking $a>0$ suitably we may suppose that $0<a \leq s$ and 
\begin{equation}
\label{eq:a, 2a, emptyset}
   \bigl\{ a, 2a, 3a, \ldots \bigr\} \cap 
   \bigl\{ \mathrm{Re} \lambda_1(0),\ldots, 
           \mathrm{Re} \lambda_m(0) \bigr\} 
   = \emptyset   
\end{equation}
hold.  By Proposition~\ref{prop:a, 2a}, 
we have such an $R_1>0$   that
the properties $\mathrm{(*)}_k$ and $\mathrm{(\sharp)}_k$ 
are valid
for $k=1,2, \ldots$.  Since $R_1>0$ can be  very small, 
we may assume that 
$a(t,x)$, $b_{j,\alpha}(t,x)$  and 
$g_{\nu}(t,x)$  have the following properties:
\begin{itemize}
\item[i)] \(
             a(t,x) \in \widetilde{\mathcal O}_a
            (S(\varepsilon (y) \times D_{R_1})
          \), 
\item[ii)] 
  \( 
     b_{j,\alpha}(t,x) \in \widetilde{\mathcal O}_a
     (S(\varepsilon (y) \times D_{R_1})
  \) 
    for $(j,\alpha) \in I_m$, $|\alpha|>0$, 
\item[iii)] 
  \( 
     g_{\nu}(t,x) \in \widetilde{\mathcal O}_0
    (S(\varepsilon (y) \times D_{R_1})
  \) 
  for 
  $|\nu| \geq 2$.
\end{itemize}

We shall construct a formal solution 
of (\ref{eq:nonlin_fuchs2}) in the form
\begin{equation}
     \label{eq:formal_sol}
     u(t,x)= \sum_{k \geq 1} u_k(t,x), 
     \quad
     u_k(t,x)
     \in \widetilde{\mathcal O}_{ak}
       \bigl(S(\varepsilon (y) \times D_{R_1})\bigr).
\end{equation}

Let us decompose our equation (\ref{eq:nonlin_fuchs2}). 
We have formally
\begin{align}
\label{eq:recurrence}
    \sum_{k \geq 1} 
    C (t \partial_t, x  )u_k 
     &=a(t,x)+ \sum_{k \geq 1} 
           b_{j,\alpha}(t,x) 
           (t \partial_t)^j\partial_x^{\alpha} u_k
       \\
     &\qquad 
     + \sum_{|\nu| \geq 2} g_{\nu}(t,x) 
      \prod_{(j,\alpha) \in I_{m}} 
      \left[ 
      \sum_{k \geq 1} (t \partial_t)^j\partial_x^{\alpha} u_k 
      \right]^{\nu_{j,\alpha}}.
      \nonumber
\end{align}
Therefore, the   equation (\ref{eq:nonlin_fuchs2}) 
is satisfied if  $\{u_k(t,x) \,;\, k=1,2,\ldots \}$ is 
determined by the following recurrent family of equations 
(\ref{eq:u_1}) and (\ref{eq:u_k}): 
\begin{equation}
\label{eq:u_1}
    C  (t \partial_t, x  )u_1
           =a(t,x) 
\end{equation}
and for $k \geq 2$, 
\begin{align}
\label{eq:u_k}
    C  (t \partial_t, x  )u_k
    =& \sum_{(j,\alpha) \in I_m} b_{j,\alpha}(t,x)
        (t \partial_t)^j\partial_x^{\alpha} u_{k-1}  
        \\
    &+ \sum_{2 \leq |\nu| \leq k} g_{\nu}(t,x)
     \sum_{|k(\nu)| = k} \prod_{(j,\alpha) \in I_m}
      \prod_{l=1}^{\nu_{j,\alpha}}
       (t \partial_t)^j\partial_x^{\alpha}u_{k_{j,\alpha}(l)}
      ,
     \nonumber
\end{align}
where 
\begin{align*}
     & k_{j, \alpha}(l) \in \mathbb{N}^* , 
     \\
    & k(\nu) = \{ (k_{j,\alpha}(l)) \,;\, (j,\alpha) \in I_m, 
           1 \leq l \leq \nu_{j,\alpha} \},
     \\
    & |k(\nu)| = \sum_{(j,\alpha) \in I_m}(k_{j,\alpha}(1)+\cdots 
     + k_{j,\alpha}(\nu_{j,\alpha})).
\end{align*}
It should be remarked that in the right hand side of 
$(\ref{eq:u_k})$ only the terms $u_1, \ldots, u_{k-1}$ and 
their derivatives appear. 
Thus, by applying  Proposition~\ref{prop:a, 2a} 
   to (\ref{eq:u_1}) and 
$(\ref{eq:u_k})$ ($k \geq 2$) 
inductively on $k$ we can obtain a 
solution $\{u_k(t,x) \,;\, k=1,2,\ldots \}$ of the recurrent 
family (\ref{eq:u_1}) and $(\ref{eq:u_k})$ (for $k \geq 2$) such that 
\begin{equation}
\label{eq:Djalpha}
    (t \partial_t)^j\partial_x^{\alpha} u_k \in 
   \widetilde{\mathcal O}_{ak}(S(\varepsilon (y) \times D_{R_1}))
   \quad \mbox{for $(j,\alpha) \in I_m$, $k=1,2,\ldots$}.
\end{equation}
This proves
\begin{proposition}
\label{pro:under}
   In  the above situation, we can construct 
\begin{equation}
\label{eq:sum/of/u_k}
     u(t,x)=\sum_{k \geq 1} u_k(t,x)
\end{equation}
with the condition (\ref{eq:Djalpha}) 
which solves the equation 
(\ref{eq:nonlin_fuchs2}) 
formally in the   sense that  
$\{u_k(t,x) \,;\, k=1,2,\ldots \}$ satisfies  
(\ref{eq:u_1}) and 
(\ref{eq:u_k}) (for $k \geq 2$).
\end{proposition}

   Set $f_1 (t)=a(t,x)$ and for $k \geq 2$
\begin{align}
\label{eq:f_k}
    f_k (t)=& \sum_{(j,\alpha) \in I_m} b_{j,\alpha}(t,x)
        (t \partial_t)^j\partial_x^{\alpha} u_{k-1} 
        \\
    &+ \sum_{2 \leq |\nu| \leq k} g_{\nu}(t,x)
     \sum_{|k(\nu)| = k} \prod_{(j,\alpha) \in I_m}
       \prod_{l=1}^{\nu_{j,\alpha}}
       (t \partial_t)^j\partial_x^{\alpha}u_{k_{j,\alpha}(l)}
       .
     \nonumber
\end{align}
Then, by $(\sharp)_k$ of Proposition~\ref{prop:a, 2a} 
we see also
\begin{proposition}
\label{prop:additional}
 In Proposition~\ref{pro:under},   we have the 
following additional property: if $f_k$ satisfies
\begin{equation}
\label{eq:additional/f_k}
    \| f_k(t) \|_r \leq C |t|^{ak} \quad 
    \mbox{on $S_{\theta}(\varepsilon (y))$}
\end{equation}
for some $0<r<R_1$, $C>0$ and $\theta>0$, we have the estimate
\begin{equation}
\label{eq:additional/tdtu_k}
    \|  (t \partial_t )^{j} 
            u_k(t) \|_r \leq 
    \frac{M_{\theta}}{k^{m-j}} \,C |t|^{ak} \enskip
    \mbox{on $S_{\theta}(\varepsilon (y))$}
    \enskip \mbox{for $j=0,1, \ldots, m-1$}.
\end{equation}
\end{proposition}

 Next,  let us prove the convergence of 
 the formal solution (\ref{eq:sum/of/u_k}).
Our aim is to show
that (\ref{eq:sum/of/u_k}) gives an $\widetilde{\mathcal O}_+$-solution of
(\ref{eq:nonlin_fuchs2}).  
\par
   Take any $\theta >0$ and $0<R<R_1$ (with $0<R \leq 1$) and fix 
them.  To show our aim it is sufficient to prove that
the formal solution (\ref{eq:sum/of/u_k}) is convergent in 
$\widetilde{\mathcal O}_a (S_{\theta}(\delta) \times D_{R/2})$
for some $\delta>0$.  
\par
   By the assumption there exist constants 
$B_{j,\alpha} \geq 0$ ($(j,\alpha) \in I_m$) and 
$G_{\nu} \geq 0$ ($|\nu| \geq 2$) satisfying the following 
properties:
\begin{itemize}
  \item[i)] $|b_{j,\alpha}(t,x)| \leq B_{j,\alpha} |t|^a$  
      on $S_{\theta}(\varepsilon (y)) \times D_R$,  
  \item[ii)]    $|g_{\nu}(t,x)| \leq G_{\nu}$ 
  on $S_{\theta}(\varepsilon (y)) \times D_R$, 
  \item[iii)] $\sum_{|\nu| \geq 2} G_{\nu} Z^{\nu}$ 
         is convergent in a neighborhood of $Z=0 \in \mathbb{C}^N$.
\end{itemize}
   By (\ref{eq:Djalpha}) we have 
$(t \partial_t)^j\partial_x^{\alpha}u_1 \in \widetilde{\mathcal O}_{a}(S(\varepsilon (y) 
\times D_{R_1}))$ for any $(j,\alpha) \in I_m$; therefore we can
take a constant $A_1 \geq 0$ such that
\begin{equation}
\label{eq:Djalpha/u_1}
   \| (t \partial_t)^j\partial_x^{\alpha}u_1(t) \|_R \leq A_1 |t|^a
   \enskip \mbox{on $S_{\theta}(\varepsilon (y))$}, \quad
   (j,\alpha) \in I_m. 
\end{equation}
Set $\beta =(e m)^m$. Using these $A_1$, $\beta$, 
$B_{j,\alpha}$   and 
$G_{\nu}$,  we consider the following holomorphic 
functional equation with respect to $Y$:
\begin{equation}
\label{eq:analytic}
    Y= A_1z + 
   M_{\theta} \left[
         \sum_{(j,\alpha) \in I_m} 
          \frac{B_{j,\alpha} z \bigl( \beta Y \bigr) }{(R-r)^m}
     + \sum_{|\nu| \geq 2}
       \dfrac{G_{\nu}}{(R-r)^{m(|\nu|-1)}}  
       \bigl( \beta Y \bigr)^{|\nu|} \right]
       ,
\end{equation}
where $r$ is a parameter with $0<r<R$.
\par
   By the implicit function theorem we see that the equation 
(\ref{eq:analytic}) has a unique holomorphic solution $Y(z)$ with $Y(0)=0$
in a neighborhood of $z=0$.  If we expand this into
\[
      Y(z)=\sum_{k \geq 1} Y_k z^k, 
\] 
we see that the coefficients
$Y_k$ ($k=1,2, \ldots$) are determined uniquely by the following
recurrence formulas:
\begin{equation}
\label{eq:Y_1}
    Y_1 = A_1  
\end{equation}
and for $k \geq 2$, 
\begin{align}
\label{eq:Y_k}
    Y_k 
    &=  M_{\theta} 
    \left[
         \sum_{(j,\alpha) \in I_m} 
          \frac{B_{j,\alpha} \beta Y_{k-1}}{(R-r)^m} \right.
    \\
    & \hspace{5em}
    + \sum_{2 \leq |\nu| \leq k}
       \frac{G_{\nu}}{(R-r)^{m(|\nu|-1)}} \left.
       \sum_{|k(\nu)| = k} 
       \prod_{(j,\alpha) \in I_m}
       \prod_{l=1}^{\nu_{j,\alpha}} 
       \beta Y_{k_{j,\alpha}(l)}
       \right] . 
       \nonumber  
\end{align}
Moreover, by induction on $k$ we   see that each $Y_k$
has the form
\begin{equation}
\label{eq:Y_k/redux}
      Y_k = \dfrac{C_k}{(R-r)^{m(k-1)}}, 
      \quad k=1,2,\ldots,  
\end{equation}
where $C_1=A_1$ and $C_k \geq 0$ ($k \geq 2$) are constants 
independent of the parameter $r$.
\par
   The following lemma guarantees that $Y(z)$ can be used as a 
majorant series of the formal solution (\ref{eq:sum/of/u_k}). 
\begin{proposition}
\label{prop:Djalpha/u_k/estimate}
For any $k=1,2, \ldots$ we have
\begin{equation}
\label{eq:Djalpha/u_k/estimate}
    \bigl\| (t \partial_t)^j\partial_x^{\alpha} u_k(t) \bigr\|_r 
       \leq \beta Y_k(r) |t|^{ak}
    \quad \mbox{on } S_{\theta}(\varepsilon (y))  
\end{equation}
for any   $0<r<R$  and $(j,\alpha) \in I_m$.
\end{proposition}

\begin{proof}
    By the definition of
$A_1$ in (\ref{eq:Djalpha/u_1}) we have
\[
     \bigl\| (t \partial_t)^j\partial_x^{\alpha} u_1(t) \bigr\|_r
     \leq A_1 |t|^a =Y_1 |t|^a\leq \beta Y_1 |t|^a
    \quad \mbox{on $S_{\theta}(\varepsilon (y))$}.
\]
This proves (\ref{eq:Djalpha/u_k/estimate}) for $k=1$. 
\par
   Let us show the general case by induction on $k$.
Suppose that $k \geq 2$ and that 
(\ref{eq:Djalpha/u_k/estimate}) has already been proved
for $u_1, \dots, u_{k-1}$. Then, by (\ref{eq:f_k}) we have
\begin{align*}
   \| f_k(t) \|_r 
    &\leq \sum_{(j,\alpha) \in I_m} 
          B_{j,\alpha} |t|^a \times 
             \beta Y_{k-1} |t|^{a(k-1)} \\
    &\qquad + 
       \sum_{2 \leq |\nu| \leq k} G_{\nu} 
       \sum_{|k(\nu)| = k} 
       \prod_{(j,\alpha) \in I_m}
       \prod_{l=1}^{\nu_{j, \alpha}}
                \beta Y_{k_{j,\alpha}(l)} 
               |t|^{a k_{j,\alpha}(l)} 
               \\
   &= |t|^{ak}
      \Biggl[ 
        \sum_{(j,\alpha) \in I_m} B_{j,\alpha} \beta Y_{k-1}  
       + 
       \sum_{2 \leq |\nu| \leq k} G_{\nu} 
       \sum_{|k(\nu)| = k} 
       \prod_{(j,\alpha) \in I_m}
       \prod_{l=1}^{\nu_{j, \alpha}}
        \beta Y_{k_{j,\alpha}(l)} 
      \Biggr] .
\end{align*}
Therefore, by comparing this with 
(\ref{eq:Y_k}) and by
using $1/(R-r)>1$ we have
\[
    \| f_k(t) \|_r \leq \frac{(R-r)^m}{M_{\theta}} Y_k |t|^{ak}
    = \frac{C_k}{M_{\theta} (R-r)^{m(k-2)}}
      |t|^{ak}
    \quad \mbox{on }  S_{\theta}(\varepsilon (y)) 
\]
for any $0<r<R$.  Hence, by Proposition~\ref{prop:additional} 
we have
\begin{equation}
  \label{eq:u_k/estimate}
     \left\| (t \partial_t)^{j} 
            u_k(t) \right\|_r \leq 
    \frac{1}{k^{m-j}} \,\dfrac{C_k}{(R-r)^{m(k-2)}}
    |t|^{ak} \enskip
    \mbox{on $S_{\theta}(\varepsilon (y))$} 
\end{equation}
for any $0<r<R$ and $j=0,1, \ldots, m-1$. 
By applying Lemma~\ref{lem:nagumo}
 (Nagumo's lemma) to this estimate we have
\begin{align*}
    \bigl\| (t \partial_t)^j\partial_x^{\alpha}u_k \bigr\|_r
    &\leq \dfrac{1}{k^{m-j}}
        \dfrac{\{m(k-2)+1\} \cdots \{m(k-2)+|\alpha|\}
        e^{|\alpha|}C_k}{(R-r)^{m(k-2)+|\alpha|}} |t|^{ak} \\
    &\leq \dfrac{1}{k^{m-j-|\alpha|}}
        \dfrac{m^{|\alpha|} e^{|\alpha|}C_k}
          {(R-r)^{m(k-2)+|\alpha|}} |t|^{ak}\\
    &\leq \dfrac{\beta C_k}
          {(R-r)^{m(k-2)+|\alpha|}} |t|^{ak} 
    \leq \dfrac{\beta C_k}
          {(R-r)^{m(k-1)}} |t|^{ak}= \beta Y_k |t|^{ak} \quad 
\end{align*}
on $S_{\theta}(\varepsilon (y))$ for any $0<r<R$ and 
$(j,\alpha) \in I_m$; this proves 
(\ref{eq:Djalpha/u_k/estimate}). 
\par
   Thus, we have proved 
   Proposition~\ref{prop:Djalpha/u_k/estimate}. 
\end{proof}

\vspace{5mm}
   Lastly, let us complete the  proof of Theorem~\ref{thm:existence_tahara}.  
Set $r=R/2$ and fix it.  Since 
$Y(z)=\sum_{k \geq 1}Y_k(r)z^k$ is convergent, we can take
a small constant $\delta >0$ so that 
$C = \sum_{k \geq 1}\beta Y_k(r)\delta^{ak} < \infty$ holds.
Then, for any $(t,x) \in S_{\theta}(\delta) \times D_{R/2}$
we have
\[
    \sum_{k \geq 1}|u_k(t,x)|
    \leq \sum_{k \geq 1}\|u_k(t) \|_{R/2}
    \leq \sum_{k \geq 1} \beta Y_k |t|^{ak}
    \leq \sum_{k \geq 1} \beta Y_k  \delta^{ak}
         \frac{|t|^a}{\delta^a}
    \leq \frac{C |t|^a}{\delta^a} .
\]
This proves that the formal solution (\ref{eq:sum/of/u_k}) is convergent in 
$\widetilde{\mathcal O}_a (S_{\theta}(\delta) \times D_{R/2})$.
\par
   Thus, we have proved Theorem~\ref{thm:existence_tahara}.

%
\end{document}